\documentclass{amsart}
\usepackage[english]{babel}
\usepackage{latexsym}
\usepackage{amssymb}
\usepackage{amsmath}
\usepackage{amscd}
\usepackage[matrix,arrow,curve]{xy}
\CompileMatrices

\DeclareMathOperator{\repq}{Rep} \DeclareMathOperator{\HOMM}{Hom}
 
 \DeclareMathOperator{\tr}{tr}
\DeclareMathOperator{\matn}{Mat} \textwidth = 450pt \textheight =
640pt \oddsidemargin = 10pt \evensidemargin = 10pt \voffset = -50 pt
\newtheorem{theor}{$\phantom{sss}$Theorem}
\newtheorem{utver}{$\phantom{sss}$Proposition}
\newtheorem{laemma}{$\phantom{sss}$Lemma}

\begin{document}
\author[S. Fedotov]{Stanislav Fedotov}
\title[Semi-invariants of 2-representations of quivers]{Semi-invariants of 2-representations of quivers}
\thanks{Supported by grant RFFI 09-01-90416 - Ukr-f-a}
\address{117463, Novoyasenevsky pr-t, 32-1-559, Moscow, Russia\endgraf
Moscow State University, department of Higher Algebra\endgraf Email:
glwrath@yandex.ru} \sloppy

\maketitle

\begin{abstract} In this work we obtain a version of the
Procesi-Rasmyslov Theorem for the algebra of semi-invariants of
representations of an arbitrary quiver with dimension vector
$(2,2,\ldots,2)$.
\end{abstract}

\section{Introduction}

We work over a base field $\Bbbk$ of characteristic zero. A quiver
$Q$ is a directed graph, determined by two finite sets $Q_0$ (the
set of ``vertices'') and $Q_1$ (the set of ``arrows'') with two maps
$h,t: Q_1\rightarrow Q_0$ which indicate the vertices at the head
and tail of each arrow. A representation $(W, \varphi)$ of $Q$
consists of a collection of finite dimensional $\Bbbk$-vector spaces
$W_v$, for each $v\in Q_0$, together with linear maps $\varphi_a:
W_{ta}\rightarrow W_{ha}$, for each $a\in Q_1$. The dimension vector
$\alpha\in\mathbb{Z}^{Q_0}$ of such a representation is given by
$\alpha_v = \dim_{\Bbbk}{W_v}$. A morphism $f:(W_v,
\varphi_a)\rightarrow (U_v,\psi_a)$ of representations consists of
linear maps $f_v: W_v\rightarrow U_v$, for each $v\in Q_0$, such
that $f_{ha}\varphi_a = \psi_af_{ta}$, for each $a\in Q_1$.
Evidently, it is an isomorphism if and only if each $f_v$ is.

Having chosen vector spaces $W_v$ of dimension $\alpha_v$, the
isomorphism classes of representations of $Q$ with dimension vector
$\alpha$ are in natural one-to-one correspondence with the orbits of
the group
$$GL(\alpha) := \prod_{v\in Q_0}GL(W_v)$$in
the representation space
$$\repq(Q,\alpha) := \bigoplus_{a\in Q_1}\HOMM(W_{ta},W_{ha}).$$

This action is given by
$(g\cdot\varphi)_a~=~g_{ha}\varphi_ag^{-1}_{ta}$, where $g =
(g_v)_{v\in Q_0}\in GL(\alpha)$. Note that the one-parameter
subgroup $\Delta = \{(tE,\ldots,tE)\}$ acts trivially.

One can also consider the action of the smaller group $SL(\alpha) =
\prod_{v\in Q_0}SL(W_v)$ on $\repq(Q,\alpha)$;
$SL(\alpha)$-invariant functions on $\repq(Q,\alpha)$ are usually
called semi-invariants. In this work we study generators of the
algebra $\Bbbk[\repq(Q,\alpha)]^{SL(\alpha)}$. Recall that
generators of $\Bbbk[\repq(Q,\alpha)]^{GL(\alpha)}$ are given by the
Procesi-Razmyslov Theorem.

\medskip

\begin{theor}~\cite[Theorem 3.4]{Claudio},~\cite{Razz} For an arbitrary
quiver $Q$ and dimension vector $\alpha$ the algebra
$\Bbbk[\repq(Q,\alpha)]^{GL(\alpha)}$ is generated by the traces of
oriented cycles of length not greater than $(\sum_v{\alpha_v})^2$.
All relations among them can be deduced from Cayley-Hamilton
polynomials.
\end{theor}

\medskip

As for the algebra $\Bbbk[\repq(Q,\alpha)]^{SL(\alpha)}$, we only
have several descriptions of its spanning sets; the main approaches
are presented in~\cite{Domokos}, \cite{DW} and~\cite{Sch}. In this
work we describe a generating set for the algebra of
2-representations of an arbitrary quiver, i.e., of representations
with dimension vector $(2,2,\ldots,2)$, in the spirit of the
Procesi-Rasmyslov Theorem.

For a matrix $A$ we will, as usually, denote by $\widehat{A}$ its
adjoint matrix, i.e., the matrix consisting of cofactors to the
elements of $A^T$. If $A_F$ is the matrix of a linear map $F:
U\rightarrow V$, then it is convenient to assume that its adjoint
matrix defines a linear map from $V$ to $U$.

Consider a quiver $Q$. Let $Q_0 = \{1,\ldots,n\}$ and $Q_1 =
\{a_1,\ldots,a_s\}$. To $Q$ we associate the quiver $\widetilde{Q}$
with $\widetilde{Q}_0 = Q_0$ and $\widetilde{Q}_1 =
\{a_1,\ldots,a_s\}\cup\{b_1,\ldots,b_s\}$, where $hb_i = ta_i$,
$tb_i = ha_i$. For each representation $(W, \varphi)$ of $Q$,
consider the {\it associated representation} of $\widetilde{Q}$ with
the same spaces $W_v$ and maps $\varphi_{a_i}$ and $\varphi_{b_i} =
\widehat{\varphi}_{a_i}$.
By a {\it route} in $Q$ we mean an oriented cycle in
$\widetilde{Q}$. For example, any cycle in $Q$ is a route. We say
that a route is {\it simple}, if no edge appears twice in the
corresponding oriented cycle. The {\it trace of a route} is the
trace of the corresponding cycle in $\widetilde{Q}$ in the
associated representation.

\begin{theor} For a quiver $Q$ and a dimension vector
$\alpha = (2, 2, \ldots, 2)$ the algebra $\Bbbk[\repq(Q,
\alpha)]^{SL(\alpha)}$ is generated by the traces of simple routes.
\end{theor}

{\bf Example 1.} Let $Q$ be the quiver $\xymatrix{
1\ar@<1ex>[r]^a\ar@<-1ex>[r]_b & 2}.$ Denote by $|X|$ the
determinant of a matrix $X$. By Theorem 2 the algebra
$\Bbbk[\repq(Q, (2,2))]^{SL(2)\times SL(2)}$ is generated by
$\tr{\varphi_a\widehat{\varphi}_a} = 2|\varphi_a|$,
$\tr{\varphi_b\widehat{\varphi}_b} = 2|\varphi_b|$ and
$\tr{\varphi_a\widehat{\varphi}_b}$ (note that
$\tr{\varphi_b\widehat{\varphi}_a} =
\tr\widehat{\varphi_a\widehat{\varphi}_b} =
\tr\varphi_a\widehat{\varphi}_b$). So, if $\varphi_a =
\left(\begin{smallmatrix}x_{11} & x_{12}\\ x_{21} & x_{22}
\end{smallmatrix}\right)$, $\varphi_b =
\left(\begin{smallmatrix}y_{11} & y_{12}\\ y_{21} & y_{22}
\end{smallmatrix}\right)$, then $$\Bbbk[\repq(Q, (2,2))]^{SL(2)\times SL(2)} =
\Bbbk[2x_{11}x_{22} - 2x_{12}x_{21}, 2y_{11}y_{22} - 2y_{12}y_{21},
x_{11}y_{22} - x_{12}y_{21} - x_{21}y_{12} + x_{22}y_{11}].$$

\medskip

Our arguments use the description of a spanning set of the algebra
$\Bbbk[\repq(Q), \alpha]^{SL(\alpha)}$ devised by Domokos and
Zubkov~\cite{Domokos}; we briefly recall it in Section 2. In
Sections 3 and 4 we obtain a formula expressing the determinant of a
$2$-block matrix as a polynomial in the traces of associated routes.
In Section 5 we prove Theorem 2.

The author hopes to use the information received about the basic
semi-invariants to generalize King's construction~\cite{Kinge} in
the way similar to~\cite{Cox}.

\smallskip

The author thanks I.V. Arzhantsev for the idea of the work and
useful discussions.

\section{The Domokos-Zubkov Theorem}

Let us recall the results of~\cite{Domokos}. Let $Q_0 =
\{1,\ldots,n\}$.
Fix a dimension vector $\alpha$ and two tuples $\overline{\imath} =
(i_1,\ldots,i_k)$ and $\overline{\jmath} = (j_1,\ldots,j_l)$ of
integers from $1$ to $n$ (possibly repeating) such that
$(\alpha_{i_1} + \ldots + \alpha_{i_k}) = (\alpha_{j_1} + \ldots +
\alpha_{j_l})$, and consider all possible matrices of size
$(\alpha_{i_1} + \ldots + \alpha_{i_k})\times (\alpha_{j_1} + \ldots
+ \alpha_{j_l})$ of form
\begin{equation*}
A_{\overline{\imath}\,\overline{\jmath}} \quad:=\quad
\begin{matrix}
 & \begin{matrix}
 i_1\phantom{A_1} & \ldots & \phantom{A_1}i_s\phantom{A_1} & \ldots & \phantom{A_1}i_k
 \end{matrix}\\
 \begin{matrix}
 j_1\\
 \vdots\\
 j_r\\
 \vdots\\
 j_l
 \end{matrix} &
\begin{pmatrix}
y_{11}F_{11} & \ldots & y_{1s}F_{1s} & \ldots & y_{1k}F_{1k}\\
\vdots & & \vdots & & \vdots\\ y_{r1}F_{r1} & \ldots & y_{rs}F_{rs}
& \ldots & y_{rk}F_{rk}\\ \vdots &  & \vdots &  & \vdots\\
y_{l1}F_{k1} & \ldots & y_{ks}F_{ls} & \ldots & y_{lk}F_{lk}
\end{pmatrix}
\end{matrix}, \eqno{(3)}
\end{equation*}
where $y_{rs}$ are formal variables, and each matrix $F_{rs}$ may be
either $0$, or the matrix of a map, that corresponds to an arrow
going from the $i_s$-th vertex to the $j_r$-th one, or an identity
matrix if $i_s = j_r$. For a fixed representation $(W, \varphi)$ a
matrix of form (3) defines a map from $W_{i_1}\oplus\ldots\oplus
W_{i_k}$ to $W_{j_1}\oplus\ldots\oplus W_{j_k}$. Moreover, its
determinant is a polynomial in variables $y_{rs}$ with
$SL(\alpha)$-invariant coefficients:
$|A_{\overline{\imath}\overline{\jmath}}| =
\sum_{\overline{\mu}}y^{\overline{\mu}}h_{\overline{\mu}}(x_{pq}^{rs})$,
where $\overline{\mu} = (\mu_{ij})_{i,j = 1}^{k, l}$ are
multidegrees, and $x^{rs}_{pq}$ are matrix elements of the matrices
$F_{rs}$.

\begin{theor}~\cite[Thm. 4.1]{Domokos} The algebra
$\Bbbk[\repq({Q},\alpha)]^{SL(\alpha)}$ is spanned by
semi-invariants $h_{\overline{\mu}}(x_{pq}^{rs})$.
\end{theor}

Thus all basic semi-invariants of ${Q}$ are given by coefficients of
monomials in variables $y_{rs}$ in determinants of block matrices of
form (3). We will describe them precisely for $\alpha =
(2,2,\ldots,2)$.

\section{Block matrices and associated routes}

For a matrix $X = \left(\begin{smallmatrix}x_{11} & x_{12}\\x_{21} &
x_{22}\end{smallmatrix}\right)$ its adjoint matrix $\widehat{X}$
equals $\left(\begin{smallmatrix}x_{22} & -x_{12}\\-x_{21} &
x_{11}\end{smallmatrix}\right)$. Recall that for $X,Y\in
\matn_{2\times 2}(\Bbbk)$ we have $\tr{X\widehat{Y}} =
\tr{\widehat{X}Y} = \tr{X}\tr{Y} - \tr{XY}$, $\tr{\widehat{X}} =
\tr{X}$ and $\widehat{XY} = \widehat{Y}\widehat{X}$.

\begin{laemma} Let $X,Y,Z\in\matn_{2\times 2}(\Bbbk)$. Then
\begin{enumerate}
\item
$\tr{X^2Y} = \tr{X}\tr{XY} - |X|\tr{Y}$;
\item
$\tr{XYXZ} = \tr{XY}\tr{XZ} - |X|\tr{Y\widehat{Z}}$.
\end{enumerate}
\end{laemma}

\begin{proof} It suffices to prove a polynomial
identity for the elements of a Zariski open subset. Hence, we may
assume that all matrices are invertible. By the Cayley-Hamilton
Theorem for two by two matrices we have
$$X^2 - (\tr{X})X + |X|E = 0.$$
Multiplying this equality by $Y$ and taking the trace of the product
received, we get (1). If $Y$ is invertible, it follows that
$$\tr{XYXZ} = \tr{((XY)^2Y^{-1}Z)} = \tr{XY}\tr{(XY\cdot Y^{-1}Z)} -
|XY|\tr{Y^{-1}Z}=$$ $\phantom{AAAAAAAAAAAAAAAAAAAAA}=
\tr{XY}\tr{(XZ)} - |X|\tr{\widehat{Y}Z}.$
\end{proof}

Consider a matrix $Z \in \matn_{2k\times 2k}(\Bbbk)$ divided into
blocks $X_{ij}$ of size $2\times 2$. To $Z$ we associate a quiver
$\Gamma$ with $\Gamma_0 = \{1,\ldots,k\}\cup\{-1,\ldots,-k\}$,
$\Gamma_1 = \{a_{ij}\mid i,j = 1,\ldots,k\}$, $ha_{ij} = i$ and
$ta_{ij} = - j$. A route in $\Gamma$ defines a sequence of matrices:
if an arrow goes from $-j$ to $i$, we take $X_{ij}$, and if it goes
from $k$ to $-l$, we take $\widehat{X}_{kl}$. Thus constructed
sequences we call {\it routes in} $Z$. The {\it adjoint} for a given
route $P = (X_1,\ldots,X_N)$ is $\widehat{P} =
(\widehat{X}_N,\ldots,\widehat{X}_1)$. Observe that any cyclic
permutation of factors as well as taking the adjoint route does not
change the trace of the product of matrices along a route. We say
that two routes are {\it equivalent} if there is such a
transformation turning one of them into another. We claim that the
determinant of $Z$ is a polynomial in traces of its routes. To prove
this we need to introduce the following construction.


\medskip

{\bf Construction.} The determinant of $Z$ equals $|Z| =
\sum_{\sigma\in S_{2k}}(-1)^{\sigma}z_{1\,\sigma(1)}\cdot
z_{2\,\sigma(2)}\cdot\ldots\cdot z_{2k\,\sigma(2k)}$. To each
summand $z_{\sigma} := z_{1\,\sigma(1)}\cdot
z_{2\,\sigma(2)}\cdot\ldots\cdot z_{2k\,\sigma(2k)}$ we assign the
{\it associated route set} $\mathcal{P}_{\sigma}$ constructed as follows.
Let $z_{1\,\sigma(1)}$ be in a block $X_{1s_1}$. 

1a) if $z_{2\,\sigma(2)}$ is in the same block as
$z_{1\,\sigma(1)}$, then add to $\mathcal{P}_{\sigma}$ the route
$(\widehat{X}_{1s_1},X_{1s_1})$ and consider the submatrix
$Z^{1\,2}_{\sigma{(1)}\,\sigma(2)}$ accessory to $X_{1s_1}$;

1b) let $z_{2\,\sigma(2)}$ be in some other $X_{1s_2}$.
In this case we consider the factor $z_{t_3\,\sigma(t_3)}$ that lies
in the same block column as $z_{2\,\sigma(2)}$, i.e., in some
$X_{r_3s_2}$; 

2) consider the factor $z_{t_4\,\sigma(t_4)}$ lying in the same
block row as the previous one, i.e., in some
$X_{r_3s_4}$;

3) continue this process as long as it is possible, alternating
horizontal and vertical shifts;

4) there comes a moment when we can not make another, some $N$-th
shift. If it were the time to move vertically (respectively
horizontally), then we have arrived to some $p$-th block column
(respectively to a $p$-th block row) for the second time. But in
each block column (except the $s_1$-th) and in each block row we
have already chosen two blocks (and hence, two factors of
$z_{\sigma}$), so it is impossible to get there again. Therefore
just before the algorithm failed we had arrived to the $s_1$-th
block column. It follows that $N$ is even (because the last shift
had been horizontal) and moreover, the last considered block was
some $X_{r_Ns_1}$;

5) consider the route $P =
(\widehat{X}_{1s_1},X_{1s_2},\widehat{X}_{r_3s_2},X_{r_3s_4},\ldots,
X_{r_ns_1})$ (all odd factors are adjoint). We add it to
$\mathcal{P}_{\sigma}$ and pass to the submatrix matrix
$Z_{\overline{P}}$ of size $(k - N)\times(k - N)$ containing all the
factors of $z_{\sigma}$ that we have not used yet.

\medskip

Note that the construction is not unique: we may take various
starting elements and choose one of two possible directions of
circuit. These transformations correspond to cyclic permutations of
the routes and/or taking the adjoint routes, hence not changing
their equivalence classes.

\medskip

{\bf Example 2.} We illustrate the construction for the block matrix
$Z = (X_{rs})_{r,s = 1}^3$ and the permutation $\sigma =
\left(\begin{smallmatrix}1 & 2 & 3 & 4 & 5 & 6\\3 & 4 & 1 & 5 & 2 &
6\end{smallmatrix}\right)$:

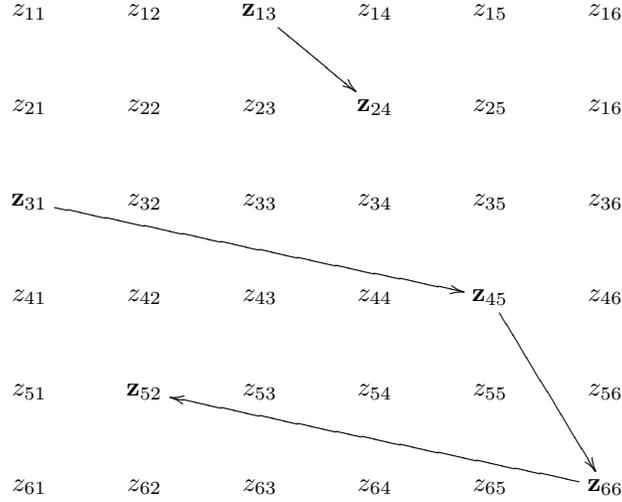
$$
\xymatrix{z_{11} & z_{12} & {\mathbf z_{13}}\ar[dr] & z_{14} & z_{15} & z_{16}\\
z_{21} & z_{22} & z_{23} & {\mathbf z_{24}} & z_{25} & z_{16}\\
{\mathbf z_{31}}\ar[drrrr] & z_{32} & z_{33} & z_{34} & z_{35} & z_{36}\\
z_{41} & z_{42} & z_{43} & z_{44} & {\mathbf z_{45}}\ar[ddr] & z_{46}\\
z_{51} & {\mathbf z_{52}} & z_{53} & z_{54} & z_{55} & z_{56}\\
z_{61} & z_{62} & z_{63} & z_{64} & z_{65} & {\mathbf
z_{66}}\ar[ullll]}
$$

The summand $z_{\sigma}$ equals
$z_{13}z_{24}z_{31}z_{45}z_{56}z_{61}$. Starting at $z_{13}$, we
next take the factor $z_{24}$. Since both of them are in the same
block $X_{12}$, we add to $\mathcal{P}_{\sigma}$ the route
$(\widehat{X}_{12},X_{12})$ and consider the submatrix
$Z^{12}_{34} = \left(\begin{smallmatrix} X_{21} & X_{23}\\
X_{31} & X_{33}\end{smallmatrix}\right)$. Here we start at $z_{31}$,
then take $z_{45}$. This element is in the second block column of
$Z^{12}_{34}$; we should take another factor lying in this column,
that is $z_{66}$, and then the factor from the same block row as
$z_{66}$; it is $z_{52}$. Thus we add to $\mathcal{P}_{\sigma}$ the
route $(\widehat{X}_{21},X_{23},\widehat{X}_{33},X_{31})$. Finally,
$\mathcal{P}_{\sigma} = \{(\widehat{X}_{12},X_{12}),\,
(\widehat{X}_{21},X_{23},\widehat{X}_{33},X_{31})\}$

\medskip

Observe that different permutations may correspond to the same
associated route sets. Namely, for $k = 1$ there are two
permutations and only one route. We say that two route sets
$\mathcal{P}$ and $\mathcal{S}$ are equivalent if their elements are
pairwise equivalent. Now, take a representative from each
equivalence class of route sets received by our construction; denote
the collection by $\mathcal{P}_k$.

\section{The determinant of a $2$-block matrix}

Define the {\it length function} $l(P)$ as the number of factors of
$P$, and the {\it index} given by
$$\nu(P)\ =\
\begin{cases} 1, &\mbox{if $P$ is of form $(\widehat{X}_{ij},X_{ij})$ or $(X_{ij},\widehat{X}_{ij})$},\\
0, &\mbox{otherwise}.
\end{cases}$$

\begin{utver} Assume that $Z = (X_{rs})_{r,s = 1}^k\in\matn_{2k\times
2k}(\Bbbk)$ is a $2$-block matrix, i.e., it is divided into blocks
$X_{rs}$ of size $2\times 2$. Then the following equality holds:
$$|Z|\ =\ \sum_{\mathcal{P}\in\mathcal{P}_k}\prod_{P\in\mathcal{P}}(-1)^{\frac12l(P) - 1}
2^{-\nu(P)}\tr{P}.\eqno{(4)}$$
\end{utver}

\begin{proof} Both sides of (4) there are polynomials in matrix elements of $X_{ij}$.
We claim that each $z_{\sigma}$ from the left side occurs in the
right side with the same coefficient.

\begin{laemma} The product $\prod_{P\in\mathcal{P}}2^{-\nu(P)}\tr{P}$
equals a sum of monomials $\pm z_{\sigma}$, where $\sigma$ run
through such permutations that every $z_{p,\sigma(p)}$ lies in an
element of some route $P\in\mathcal{P}$.
\end{laemma}

\begin{proof} First prove the lemma for $\mathcal{P} =
\{P\}$, $\nu(P) = 0$.

Note that a transposition of block rows or block columns changes
neither $|Z|$, nor the right side of (4). Hence we may assume that
$P =
(\widehat{X}_{11},X_{12},\widehat{X}_{22},X_{23},\ldots,\widehat{X}_{kk},X_{k1})$.

The polynomial $\tr{P}$ is a sum of monomials in matrix elements of
blocks $X_{ij}$ of the following form
$$(\widehat{X}_{11})_{p_1p_2}(X_{12})_{p_2p_3}(\widehat{X}_{22})_{p_3p_4}(X_{23})_{p_4p_5}
\ldots(\widehat{X}_{k,k})_{p_{2k-1}p_{2k}}(X_{k,1})_{p_{2k}p_1}.\eqno{(5)}$$

Such a product equals some $z_{\sigma}$ whenever each pair of its
factors does not lie in one row or in one column. By construction of
$P$ only some $(X_{q,q+1})_{p_{2q}p_{2q + 1}}$ and
$(\widehat{X}_{q+1,q+1})_{p_{2q+1}p_{2q+2}} =
(X_{q+1,q+1})_{\tilde{p},\tilde{q}}$ may lie in the same block
column. They are in the same column if and only if $\tilde{q} =
2q+1$. But for $X, Y\in\matn_{2\times 2}(\Bbbk)$ we have
$$X\widehat{Y} = \begin{pmatrix}
x_{11}y_{22} - x_{12}y_{21} & - x_{11}y_{12} + x_{12}y_{11}\\
x_{21}y_{22} - x_{22}y_{21} & - x_{21}y_{12} + x_{22}y_{11}
\end{pmatrix},$$
and $b\ne d$ holds in all the products $x_{ab}y_{cd}$ that occur
here. The same argument may be used to prove that every two factors
of (5) do not lie in the same row.

If $\nu(P) = 1$, i.e., $P = (\widehat{X}_{rs},X_{rs})$, we have
$\tr{P} = 2|X_{rs}|$. This gives rise to the factors $2^{-\nu(P)}$
in (4).

Assume that there is more than one route in $\mathcal{P}$. Each
route $\mathcal{P} = \{P_1,\ldots,P_m\}$ determines a collection of
block submatrices $Z_{P_1},\ldots,Z_{P_m}$ such that each $P_i$ is
entirely situated in $Z_{P_i}$ and each block row and each block
column of $Z$ intersects with a unique $Z_{P_i}$. By transpositions
of block rows and columns we can make all the blocks $Z_{P_i}$
diagonal. Now it is evident that
$\prod_{P\in\mathcal{P}}2^{-\nu(P)}\tr{P}$ equals a sum of products
$z_{\sigma_1}\ldots z_{\sigma_m}$, where $z_{\sigma_q}$ are summands
of $|Z_{P_q}|$ and $\sigma_j$ is associated to $P_j$.
\end{proof}

Thus the right side of (4) is a linear combination of $z_{\sigma}$.

\begin{laemma} Each monomial $z_{\sigma}$ occurs in
$\prod_{P\in\mathcal{P_{\sigma}}}2^{-\nu(P)}\tr{P}$ with nonzero
coefficient.
\end{laemma}

\begin{proof} Without loss of generality, we may assume that
$\mathcal{P}_{\sigma} = \{P\}$ with $P = (\widehat{X}_{r_1s_1},
X_{r_1s_2},\widehat{X}_{r_3s_2}, \ldots,\\ X_{r_{2k}s_1})$. Denote
by $z_{\sigma,q}$ the factor of $z_{\sigma}$ lying in the $q$-th
element of $P$. Fixing the row containing $z_{\sigma,1}$, we
determine the row, in which $z_{\sigma,2}$ lies: it is the second
row of a pair ($2r_1 - 1$, $2r_1$). Similarly, the choice of the
column containing $z_{\sigma,2}$ determines in which column
$z_{\sigma,3}$ lies, and so on. Finally, knowing the column
containing $z_{\sigma,2k}$ we learn in which column $z_{\sigma,1}$
lies. Consequently, all the permutations $\tau$ with
$\mathcal{P}_{\tau} = \mathcal{P}_{\sigma}$ are parametrised by
tuples $\xi(\tau) = (\xi_1,\ldots,\xi_k)\in(\mathbb{Z}_2)^{2k}$.
Further, if we interchange the $(2r_1 - 1)$-th and the $2r_1$-th
rows of $Z$, then $\xi_1$ becomes $1-\xi_1$ and the other elements
of a tuple do not change. Thus there exists a sequence of such
permutations sending a product $z_{\sigma}$ to any $z_{\tau}$ with
the same associated set. It is only left to understand, how these
transformations change $\tr{P}$. 

If we interchange the $(2s_q - 1)$-th and the $2s_q$-th columns, it
 only influences the fragment
$(X_{r_{q-2}s_q},\widehat{X}_{r_qs_q})$. But for $X, Y\in
\matn_{2\times 2}(\Bbbk)$ we have
$$X\widehat{Y} = \begin{pmatrix}
x_{11} & x_{12}\\
x_{21} & x_{22}\end{pmatrix} \begin{pmatrix} y_{22} & -y_{12}\\
-y_{21} & y_{11}\end{pmatrix} =
\begin{pmatrix}
x_{11}y_{22} - x_{12}y_{21} & - x_{11}y_{12} + x_{12}y_{11}\\
x_{21}y_{22} - x_{22}y_{21} & - x_{21}y_{12} + x_{22}y_{11}
\end{pmatrix},$$
and after the permutation:
$$\begin{pmatrix}
x_{12} & x_{11}\\
x_{22} & x_{21}\end{pmatrix} \widehat{\begin{pmatrix} y_{12} & y_{11}\\
y_{22} & y_{21}\end{pmatrix}} =
\begin{pmatrix}
x_{12}y_{21} - x_{11}y_{22} & - x_{12}y_{11} + x_{11}y_{12}\\
x_{22}y_{21} - x_{21}y_{22} & - x_{22}y_{11} + x_{21}y_{12}
\end{pmatrix}.$$
Comparing the results we see that both products contain the same
members, yet with different signs. So a transposition multiplies
$\tr{P}$ by $(-1)$. The same argument works for transformations of
columns. Hence all $z_{\tau}$ with $\mathcal{P}_{\tau} = \{P\}$
occur in $\tr{P}$ with nonzero coefficient, which in fact equals
$\pm 1$.
\end{proof}

It is clear that the same $z_{\sigma}$ can not occur as a summand in
$\prod_{S\in\mathcal{S}}2^{-\nu(S)}\tr{S}$ and
$\prod_{P\in\mathcal{P}}2^{-\nu(P)}\tr{P}$ for two different sets
$\mathcal{P}$ and $\mathcal{S}$ in $\mathcal{P}_{k}$. Therefore each
$z_{\sigma}$ occurs in both sides of (4) with coefficient $\pm 1$.
It remains to prove that these signs are the same.

Assume as before that $\mathcal{P}_{\sigma} = {P}$, where $P =
\widehat{X}_{11}X_{12}\widehat{X}_{22}\ldots X_{l(P),1}$. We know
that if $\mathcal{P}_{\tau} = \mathcal{P}_{\sigma}$, then there
exists a permutation of rows and columns of $Z$ multiplying both
sides of (4) by the same number $1$ or $(-1)$ and transforming
$z_{\sigma}$ into $z_{\tau}$. So the question is whether the
coefficients of a fixed $z_{\sigma_0}$ are the same. Consider
$\sigma_0 = (1, 2, 3, \ldots, 2l(P) - 1 , 2l(P))$. Then the
coefficient of $z_{\sigma_0}$ in $|Z|$ equals $(-1)^{\sigma_0} =
-1$. As for the right side, multiply the summand
$$\pm z_{12}z_{23}\ldots z_{2l(P) - 1,2l(P)}z_{2l(P),1} = (\widehat{X}_{11})_{21}(X_{12})_{21}(\widehat{X}_{22})_{21}(X_{23})_{21}\ldots
(\widehat{X}_{l(P),l(P)})_{21}(X_{l(P),1})_{21}$$ by $(-1)^{\frac12
l(P) - 1}$. We have $(X_{q,q+1})_{21} = z_{2q, 2q + 1}$ and
$(\widehat{X}_{qq})_{21} = - z_{2q - 1, 2q}$, so the total sign is
 $(-1)^{\frac12 l(P) - 1}(-1)^{\frac12 l(P)} =
-1$. This completes the proof of Proposition 1.
\end{proof}

\section{Proof of Theorem 2}

Recall that by Theorem 3 the algebra $\Bbbk[\repq({Q},
\alpha)]^{SL(\alpha)}$ is generated by the coefficients of monomials
in $y_{rs}$ in determinants of block matrices of form (3). So it is
only left to prove that these are precisely the traces of routes
associated to those matrices.

It is easy to see that in the determinant of
$\left(y_{rs}X_{rs}\right)_{r,s = 1}^k$ the coefficient of
$y_{r_1s_1}\ldots y_{r_ks_k}$ equals the alternating sum of all such
$z_{\sigma}$  that there exists a $\tau\in S_{2k}$, for which every
$z_{p\,\sigma(p)}$ is in $X_{r_{\tau(p)}s_{\tau(p)}}$. Observe that
each block row and each block column contains precisely two blocks
$X_{r_ps_p}$ and hence these blocks are elements of a route from
some set $\mathcal{P}$. Now from the proof of Proposition 1 it is
clear that the coefficient of $y_{r_1s_1}\ldots y_{r_ks_k}$ equals
$\prod_{P\in\mathcal{P}}(-1)^{\frac12l(P) - 1} 2^{-\nu(P)}\tr{P}$.

Now prove that traces of the routes under consideration are
semi-invariants. If $\widehat{X}_1,\ldots,\widehat{X}_N$ are all the
adjoint factors of a route $P$ (with multiplicities), then
$$\tr{P} = \tr{\widehat{X}_1P_1\ldots\widehat{X}_NP_N} =
|X_1|\cdot\ldots\cdot|X_N|\tr{X_1^{-1}P_1\ldots X_N^{-1}P_N},$$
where $P_i$ are products of the remaining elements of $P$. Since
$\tr{X_1^{-1}P_1\ldots X_N^{-1}P_N}$ is a rational invariant,
$\tr{P}$ and $|X_1|\cdot\ldots\cdot |X_N|$ have the same weight.

Further, if a factor $X$ appears twice in a route $P$, then by Lemma
1 we have
$$\tr{P} = \tr{XP_1XP_2} = \tr{XP_1}\tr{XP_2} -
\frac12\tr{\widehat{X}X}\tr{P_1\widehat{P}_2}.$$ The sequence
$(P_1,\widehat{P}_2)$ is a route yielding that
$\tr{P_1\widehat{P}_2}$ is a semi-invariant. So the algebra of
semi-invariants is generated by traces of simple routes. Theorem 2
is proved.

\medskip

Observe that the trace of a route containing a pair
$\widehat{\varphi}_a$, $\widehat{\varphi}_b$ with $ha = tb$, $ta =
hb$ may be excluded from the generating set. Indeed, such a route is
of form $(\widehat{\varphi}_a,P,\widehat{\varphi}_b,S)$, where $P$
and $S$ are subroutes. But for $X,Y\in\matn_{2\times 2}(\Bbbk)$
Lemma 1 implies that $|X|\tr{Y\widehat{Z}} = \tr{XY}\tr{XZ} -
\tr{XYXZ}$. Therefore,
$$\tr{\widehat{\varphi}_aP\widehat{\varphi}_bS} =
\frac1{|\varphi_a|}|\varphi_a|\tr{\widehat{\varphi}_aP\widehat{\widehat{S}\varphi_b}}
=
\frac1{|\varphi_a|}(\tr{\varphi_a\widehat{\varphi}_aP}\tr{\varphi_a\widehat{S}\varphi_b}
- \tr{\varphi_a\widehat{\varphi}_aP\varphi_a\widehat{S}\varphi_b})
=$$
$$= \frac1{|\varphi_a|}(|\varphi_a|\tr{P}(\tr{\varphi_a\varphi_b}\tr{S}
- \tr{\varphi_b\varphi_aS}) -
|\varphi_a|(\tr{\varphi_bP\varphi_a}\tr{S} -
\tr{\varphi_bP\varphi_aS})) =$$$$=
\tr{\varphi_a\varphi_b}\tr{P}\tr{S} - \tr{P}\tr{\varphi_b\varphi_aS}
- \tr{S}\tr{\varphi_a\varphi_bP} + \tr{\varphi_bP\varphi_aS}.$$ If
the initial product $\widehat{\varphi}_aP\widehat{\varphi}_bS$ were
assigned to a route, then $P$ and $S$ should map from $W_{ha}$ to
$W_{ha}$ and from $W_{ta}$ to $W_{ta}$ respectively:
$$
\xymatrix{W_{ta}\ar@/_1pc/[rrr]_{S} &
W_{ha}\ar[l]_{\widehat{\varphi}_a} & W_{ha}\ar[l]_{P} &
W_{ta}\ar[l]_{\widehat{\varphi}_{b}}}.$$ So the products whose
traces we take in the final formula are also associated to routes.

\medskip

Finally, note that for 3-representations Theorem 2 does not hold.
Indeed, consider the quiver

$$
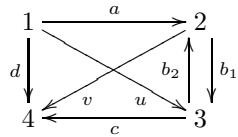
\xymatrix{1\ar[rr]^{a}\ar[d]_{d}\ar[drr]_(0.7){u} &&
2\ar@<1ex>[d]^{b_1}\ar[dll]^(0.7){v}\\
4 && 3\ar@<1ex>[u]^{b_2}\ar[ll]^{c}}
$$
and the dimension vector $(3, 3, 3, 3)$. The coefficient $F$ of
$\prod_{i,j}y_{ij}$ in the determinant of the block matrix
$$Z =
\begin{pmatrix} y_{11}\varphi_{a} & y_{12}E &
y_{13}\varphi_{b_2}\\
y_{21}\varphi_{u} & y_{22}\varphi_{b_1} & y_{23}E\\
y_{31}\varphi_{d} & y_{32}\varphi_{v} & y_{33}\varphi_{c}
\end{pmatrix}$$
is not invariant, and its total degree with respect to the variables
of each $\varphi$ equals 1. On the other hand, the trace of a route
is not invariant if and only if the route contains an adjoint
matrix. The total degree of a trace with respect to the variables
from this matrix is at least 2. Thus the subalgebra generated by the
traces of routes does not contain $F$.

\end{document}